\documentclass[a4paper,11pt]{article}

\usepackage[tbtags]{amsmath}
\usepackage{amsfonts,amssymb,amsthm}
\usepackage{eucal}
\usepackage{setspace}
\usepackage[margin=1.25in]{geometry}

\usepackage{natbib}
\usepackage{graphicx}

\usepackage{titlesec}
\titleformat{\section}{\large\bf}{\thesection}{1em}{}

\usepackage{dsfont} 
\usepackage{verbatim}
\usepackage{url}

\newcommand{\ldef}{\mathrel{:=}\nolinebreak}
\newcommand{\rdef}{\mathrel{=:}\nolinebreak}
\newcommand{\Set}[1]{\left\lbrace \boost #1 \right\rbrace}
\newcommand{\boost}{\vphantom{\big(}}
\newcommand{\bfrac}[2]{\frac{\boost #1}{\boost #2}}
\newcommand{\fbH}{f_{\bH}}
\newcommand{\fbX}{f_{\bX}}
\newcommand{\fstar}{f^*}
\newcommand{\bfstar}{\bs{f^*}}
\newcommand{\ftilde}{\tilde{f}}
\newcommand{\bfhat}{\bs{\hat{f}}}
\newcommand{\fhat}{\hat{f}}
\newcommand{\bftilde}{\bs{\tilde{f}}}
\newcommand{\given}{\mathbin{\vert}\nolinebreak}
\newcommand{\fhtilde}{{\hat{f}'}}
\newcommand{\bfhtilde}{\bs{\fhtilde}}
\newcommand{\deq}{\mathrel{\dot =}\nolinebreak}

\newcommand{\TV}[1]{\big\Vert #1 \big\Vert_{\text{TV}}}
\newcommand{\abs}[1]{\left\lvert \boost #1 \right\rvert}
\newcommand{\sumnl}{\sum\nolimits}
\newcommand{\indic}{\mathds{1}}
\newcommand{\assign}{\leftarrow}

\newcommand{\ptilde}{\tilde{p}}
\newcommand{\bptilde}{\bs{\tilde{p}}}
\newcommand{\iid}{\buildrel \rm iid \over \sim}

\newcommand{\simplex}{\mathds{S}}
\newcommand{\R}{\mathds{R}}
\newcommand{\half}{\tfrac{1}{2}}
\newcommand{\hlo}{\underline{h}}
\newcommand{\hhi}{\overline{h}} 
\newcommand{\cdott}{\cdot\,}

\newcommand{\bs}[1]{\boldsymbol{#1}}
\newcommand{\bone}{\bs{1}}

\newcommand{\bh}{\bs{h}}
\newcommand{\bH}{\bs{H}}
\newcommand{\curlyH}{\mathcal{H}}
\newcommand{\curlyM}{\mathcal{M}}
\newcommand{\bp}{\bs{p}}
\newcommand{\bphi}{\bs{\phi}}
\newcommand{\bu}{\bs{u}}
\newcommand{\curlyU}{\mathcal{U}}
\newcommand{\bw}{\bs{w}}
\newcommand{\bx}{\bs{x}}
\newcommand{\bX}{\bs{X}}
\newcommand{\bz}{\bs{z}}
\newcommand{\bZ}{\bs{Z}}
\newcommand{\curlyX}{\mathcal{X}}

\DeclareMathOperator{\E}{E}
\let\Pr\undefined
\DeclareMathOperator{\Pr}{Pr}
\newtheorem{theorem}{Theorem}

\title{Exchangeability, the `Histogram Theorem',\\and population inference}

\author{Jonathan Rougier\thanks{School of Mathematics, University Walk,
Bristol BS8 1TW, UK; email \texttt{j.c.rougier@bristol.ac.uk}.}\\School of
Mathematics\\University of Bristol}

\date{}

\begin{document}

\maketitle

\begin{abstract}\noindent%
Some practical results are derived for population inference based on a
sample, under the two qualitative conditions of `ignorability' and
exchangeability.  These are the `Histogram Theorem', for predicting the
outcome of a non-sampled member of the population, and its application
to inference about the population, both without and with groups.  There
are discussions of parametric versus non-parametric models, and different
approaches to marginalisation.  An Appendix gives a self-contained proof
of the Representation Theorem for finite exchangeable sequences.

\medskip\noindent%
\end{abstract}

\section{Introduction}
\label{sec:intro}

This paper concerns the relationship between the histogram of a sample
and proportions in the underlying population.  Our intuition is that
these two distinct objects must be similar, but cannot be identical.  As
a simple illustration, we would not automatically assign a population
proportion of zero to a label which was not present in the sample
histogram.  The purpose of this paper is to derive and publicize some
practical results: the `Histogram Theorem' (Theorem~\ref{thm:HT} in
section~\ref{sec:HT}), and its implications for sample-based population
inference, eqs~\eqref{eq:mayor} and \eqref{eq:strat}.

Any quantitative assessment of the relationship between the sample
histogram and the population must depend on a statistical model of the
sampling procedure.  I will assume that the sampling process is
`ignorable', as described in \citet[ch.~8]{gelman14}.  On the basis of
`ignorability' alone, it is possible to make inferences about the
population using `design-based' estimation, e.g.\ as discussed in
\citet{little04}, \citet{brewer02}, and \citet{brewer09}.  I will adopt
what is often characterised as a competing mode of inference, which is
to include an explicit statistical model for the population.
\citet{brewer09} term this `prediction-based' estimation, although
`model-based' is also common.

My statistical model for the population is that it is exchangeable in
the quantities of interest, either \textit{in toto} or within
identifiable groups.  A set of random quantities is exchangeable exactly
when beliefs (probabilities) about any subset are the same as beliefs
about any other subset of the same size; see section~\ref{sec:HT} for
the formal definition.  I prefer to regard exchangeability as a
`boundary condition'.  To say ``I am treating my beliefs about the
population as exchangeable'' asserts that I am choosing to ignore
everything about the population except for the values that are collected
in the sample.  Where exchangeability is confined within groups, the
assertion is that I am choosing to ignore everything about the
population except its group structure and the values that are collected
in the sample.  It is often expedient to make such an assertion; in some
circumstances it is also equitable.  On this interpretation, the
restrictions implied by the exchangeable model are to be welcomed,
rather than seen as `unrealistic'.

The key result in the paper is the `Histogram Theorem', stated and
proved in section~\ref{sec:HT}.  This result relates the predictive
probability distribution of an unsampled member of the population to the
histogram of the sample, in terms of an upper bound on the total
variation distance.  Its practical implications are discussed in
section~\ref{sec:impl}.  Population inference is described in
section~\ref{sec:pop}, which improves the folk theorem that the sample
histogram predicts the population proportions.  This result is extended
to a more limited form of exchangeability (within but not between
groups).  Section~\ref{sec:param} discusses the role of parametric
population models, and clarifies the meaning of `non-parametric' in the
context of exchangeable beliefs.  Section~\ref{sec:merge} considers the
effect of merging labels, which happens implicitly whenever a
multiple-question survey is analysed marginally, one question at a time.
A self-contained Appendix states and proves the Finite Representation
Theorem for exchangeable random quantities (Theorem~\ref{thm:FRT}),
providing more details for some of the steps in the main text.

\section{Prediction in finite exchangeable sequences}
\label{sec:HT}

Let $\bX \ldef (X_1, \dots, X_m)$ be a finite sequence of random quantities,
where each $X_i$ has the same finite realm
\begin{displaymath}
  \curlyX \ldef \Set{ x^{(1)}, \dots, x^{(k)} } .
\end{displaymath}
The finite $k$ does not exclude the case where the the realm of the
$X$'s is non-finite.  In this case the set $\curlyX$ represents a
specified finite partition of the realm of $X_i$, and may thus have
additional topological structure.  But I will not assume any such
structure in what follows, and hence I will treat the $x^{(j)}$ simply
as `labels'.

In this paper $\bX$ is an exchangeable sequence, which is to say that
\begin{equation}\label{eq:EXdef}
  \E\{ g(X_1, \dots, X_m) \} = \E\{ g(X_{\pi_1}, \dots, X_{\pi_m}) \}
\end{equation}
for any real-valued function $g$ and any permutation $(\pi_1, \dots,
\pi_m)$.  See \citet[ch.~1]{schervish95}, \citet{kingman78}, or
\citet{aldous85} for results and insights about exchangeable sequences,
and \citet{diaconis77} and \citet{diaconis80} for the special case of
finite exchangeability, as used here.  An Appendix at the end of the
paper provides the necessary details about exchangeability, including a
statement and proof of the Finite Representation Theorem
(Theorem~\ref{thm:FRT}).

An equivalent statement for exchangeable $\bX$ with finite realms is
that the probability mass function (PMF) depends only on the histogram
of $\bx$, denoted $$\bh(\bx) \ldef \big( h_1(\bx), \dots, h_k(\bx)
\big),$$ where $h_j(\bx)$ is the number of elements of $\bx$ with label
$x^{(j)}$.  Thus the PMF for $\bX$ can be expressed as
\begin{equation}\label{eq:fbXm}
  \fbX(\bx) = \fbH^m \{ \bh(\bx) \} \big/ \, \curlyM_{\bh(\bx)} 
\end{equation}
where $\fbH^m$ is a PMF for histograms of $m$ elements allocated over
$k$ labels, and $\curlyM_{\bh(\bx)}$ is the Multinomial coefficient,
defined in \eqref{eq:curlyM}, which represents the number of
distinct sequences $\bx$ which have the same histogram~$\bh(\bx)$.

If $\bX$ is an exchangeable sequence then any subsequence of $\bX$ is
also an exchangeable sequence (this is obvious from eq.~\ref{eq:EXdef}).
Without loss of generality, I will focus on the subsequence comprising
the first $n+1$ elements of $\bX$, denoted $(\bX_{1:n}, X_{n+1})$.  The
predictive distribution of $X_{n+1}$ given $\bx_{1:n}$ is denoted
\begin{equation}\label{eq:fstar0}
  \fstar_j \ldef \Pr\{ X_{n+1} \deq x^{(j)} \given \bX_{1:n} \deq \bx_{1:n} \} \qquad j = 1, \dots, k ,
\end{equation}
to be derived below; I prefer to use dots to indicate binary predicates
in infix notation.\footnote{That is, I make a notational distinction
between a statement such as `$x = 1$' which is an assertion about $x$,
and `$x \deq 1$', which is a sentence from first-order logic which is
either False or True.} This predictive distribution has an exact
expression originating from $\fbH^m$, although this expression will
typically be very complicated (see the Appendix).


One intuitive approximation to $\fstar_j$ is
\begin{equation}\label{eq:ftilde}
  \ftilde_j \ldef \frac{ h_j + 1 }{ n + k } \qquad 
  j = 1, \dots, k ,
\end{equation}
where from now on I will suppress the `$(\bx_{1:n})$' argument on
$\bh(\bx_{1:n})$ and $h_j(\bx_{1:n})$, to save clutter.  In other words,
add one to each count in the histogram and normalise the result.  In
Machine Learning this is sometimes termed `add one smoothing'
\citep[see, e.g.,][p.~79]{murphy12}.  Obviously \eqref{eq:ftilde} is a
very attractive approximation, if it is accurate, because it makes no
reference to $\fbH^m$, and is trivial to compute.  The first objective of
this paper is to prove the following result, on the basis of which I
refer to \eqref{eq:ftilde} as the `HT approximation'.  I write $\bh
\oplus j$ to denote the histogram $\bh$ with $1$ added to the count of
the $j$th label.

\begin{theorem}[Histogram Theorem]\label{thm:HT}\samepage
Let $\beta$ satisfy
\begin{displaymath}
1 + \beta = \bfrac{ \max_j \fbH^{n+1}(\bh \oplus j) }{ \min_{j'} \fbH^{n+1}(\bh \oplus j') } .
\end{displaymath}
Then the total variation distance between $(\fstar_1, \dots, \fstar_k)$
and $(\ftilde_1, \dots, \ftilde_k)$ is no greater than~$\half \, \beta$.  
\end{theorem}

Note that the Histogram Theorem is a pure exchangeability result which
requires no additional structure; e.g.\ no super-populations or
parameters, as used in \citet{ericson69}.  Parameters will be discussed
further in section~\ref{sec:param}.  An operational
interpretation of the Histogram Theorem is given at the start of
section~\ref{sec:impl}.

The proof comes in two parts.  First, I derive an exact expression for
$\fstar_j$, then I derive the total variation bound.  For the first part,
I generalise the approach of David Blackwell, in his discussion of
\citet{diaconis88}.  The marginal distribution of $\bX_{1:n}$ can be
derived from $\fbX$ given in \eqref{eq:fbXm}; denote it as
\begin{equation}\label{eq:fbXn}
  f_{\bX_{1:n}}(\bx_{1:n}) = f_{\bH}^n (\bh) \big/ \, \curlyM_{\bh} ,
\end{equation}
where $f_{\bH}^n(\bh)$ is deduced from $f_{\bH}^m$\/; this is the
necessary form because$\bX_{1:n}$ is exchangeable (see the Appendix).
The predictive distribution has its standard quotient form: 
\begin{equation}\label{eq:fstar1}
  \fstar_j = \frac{ f_{\bX_{1:n}, X_{n+1}}(\bx_{1:n}, x^{(j)}) }{ f_{\bX_{1:n}}(\bx_{1:n}) } 
\end{equation}
where we can assume that the denominator is non-zero (otherwise $\fbH^m$
would need to be revised in the light of the sample).  Write
\begin{equation}\label{eq:fstar2}
  r_j \ldef \frac{ \fstar_j }{ \fstar_1 } \quad \text{so that} \quad
  \fstar_j = \frac{ r_j }{ \sum_{j'} r_{j'} } .  
\end{equation}
Now evaluate $r_j$ in terms of \eqref{eq:fbXn} and \eqref{eq:fstar1} to
give
\begin{equation}
  r_j = \bfrac{ f_{\bH}^{n+1}(\bh \oplus j) \big/ \,
  \curlyM_{\bh \oplus j} }{ f_{\bH}^{n+1}(\bh \oplus 1) \big/ \,
  \curlyM_{\bh \oplus 1} } 
= \bfrac{ f_{\bH}^{n+1}(\bh \oplus j) }{
  f_{\bH}^{n+1}(\bh \oplus 1) } \cdot \frac{ h_j + 1 }{ h_1 + 1 } \rdef
  u_j \cdot v_j \, , \label{eq:betaj}
\end{equation}
say, introducing the terms $u_j$ and $v_j$.  In these terms,
\begin{equation}\label{eq:fstar3}
 \fstar_j = \frac{ u_j \, v_j }{ \sum_{j} u_{j'}\, v_{j'} }
 \quad\text{and}\quad \ftilde_j = \frac{ v_j }{ \sum_{j'} v_{j'}
 }
\end{equation}
from \eqref{eq:fstar2} and \eqref{eq:ftilde}, respectively.  From now on
I write 
\begin{displaymath}
  \bftilde \ldef \big( \ftilde_1, \dots, \ftilde_k \big) ,
\end{displaymath}
and similarly for $\bfstar$ above, $\bfhat$ and $\bfhtilde$ in
section~\ref{sec:impl}, and $\bptilde$ and $\bptilde_g$ in
section~\ref{sec:pop}.

The second part of the proof consists of showing under what conditions
the $u$'s in \eqref{eq:fstar3} can be ignored, allowing $\bftilde$ to
be a good approximation to $\bfstar$.  Here I co-opt the more general
result of L.J.~Savage in \citet{edwards63}, his `Principle of stable
estimation'.  In Savage's notation, my condition in the Histogram
Theorem reads
\begin{displaymath}
  \varphi \leq u_j \leq (1 + \beta) \, \varphi
\end{displaymath}
for every $u_j$, where in my case $\varphi \ldef \min_j u_j$.  It follows
immediately that
\begin{equation}\label{eq:bound}
  \varphi \sumnl_j v_j \leq \sumnl_j u_j \, v_j \leq (1 + \beta) \, \varphi \sumnl_j v_j
\end{equation}
and, with only a little more work, that
\begin{equation}\label{eq:blim}
  \frac{1}{1 + \beta} \leq \frac{ \fstar_j }{ \ftilde_j } \leq 1 + \beta \qquad j = 1, \dots, k ,
\end{equation}
from \eqref{eq:fstar3} and \eqref{eq:bound}.  Then, for the total variation distance,
\begin{align*}
  \TV{\bfstar - \bftilde}
  & \ldef \sup_{C \subset \curlyX} \abs{ \fstar(C) - \ftilde(C) } \\
  & = \frac{1}{2} \sumnl_j \abs{ \fstar_j - \ftilde_j } && \text{as $\curlyX$ is finite} \\
  & = \frac{1}{2} \sumnl_j \left\lvert \frac{ \fstar_j }{ \ftilde_j } - 1 \right\rvert \ftilde_j && \text{as $\ftilde_j > 0$} \\
  & \leq \frac{1}{2} \sumnl_j \max \left\{ \frac{\beta}{1 + \beta}, \beta \right\} \ftilde_j && \text{from \eqref{eq:blim}} \\
  & = \frac{1}{2} \, \beta 
\end{align*}
as was to be shown.  This completes the proof of the Histogram Theorem.

Savage provided a more general result than the version used here, with
conditions on the $v$'s as well as the $u$'s, that might be useful in
establishing a tighter upper bound if the histogram is highly
concentrated in a subset of the labels.

Michael Goldstein (pers.\ comm.) has provided me with a simple
illustration of when $\beta$ is infinite: the case where you are sure
that all the $X$'s take the same value, but unsure of what that value
is.  In this case the histogram $\bh$ will have all $n$ cases with the
same label, and one of the `add one' histograms adjacent to $\bh$ will
be consistent with this, while the other $k-1$ will not.  Hence $\beta$
is of the form $1 / 0$, and $\TV{\bfstar - \bftilde}$ takes its maximum
value of $1$.  But in this case a sample larger than one is not required
(except for caution), and exact predictions are straightforward.

\section{Some simple implications}
\label{sec:impl}

The Histogram Theorem (HT) suggests the following procedure for making a
prediction about an unobserved $X_i$, based on a sample taken from a
population treated as exchangeable.  First, contemplate the sample
histogram $\bh$.  If you have no strong beliefs about the `add one'
histograms adjacent to $\bh$, in the sense that you do not believe that
the most probable of them, \textit{a priori}, is much more probable than
the least probable, then your $\beta$ in the HT is small, and the HT
approximation in \eqref{eq:ftilde} gives an accurate prediction.  This
procedure is illustrated in Figure~1.

\begin{figure}[p]

\begin{center}
\begin{minipage}{0.8\textwidth}

\includegraphics[width=\textwidth]{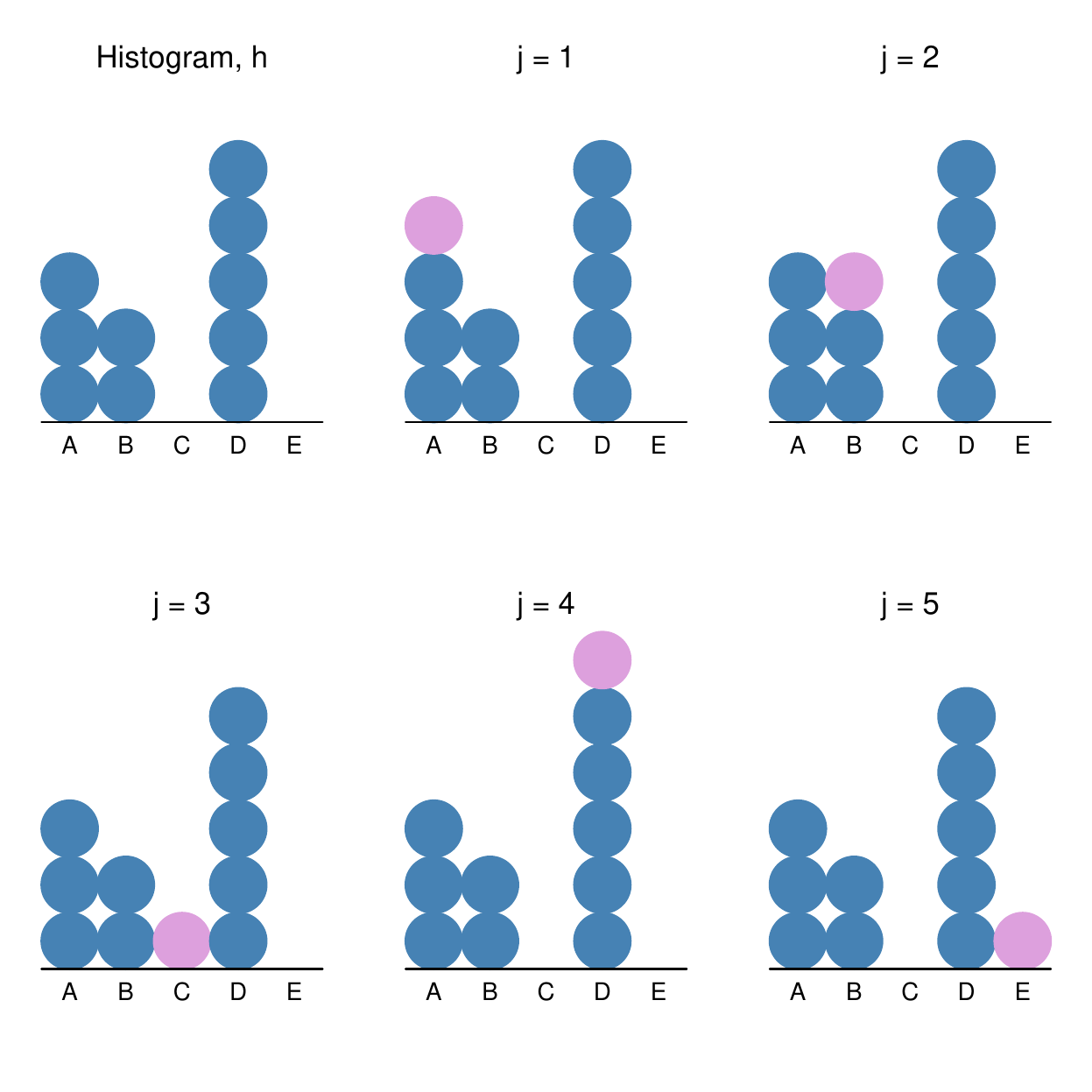}

\caption{Contemplating the histogram $\bh$ in order to reflect on the
size of $\beta$.  In this illustration, $k = 5$, $n = 10$, and $\bh =
(3, 2, 0, 5, 0)$.  The five `add one' histograms adjacent to $\bh$ are
constructed, and the value of $\beta$ is inferred from the ratio of the
probabilities of the the most probable to the least probable, \textit{a
priori}, for an exchangeable population of size $m$ (not explicitly
specified here).}

\end{minipage}
\end{center}

\label{fig:allHist}

\end{figure}

You may be comfortable with the qualitative assessment of $\beta$ as
`small'.  But if, after contemplation of $\bh$, you are prepared to
quantify $\beta$ as, say, `not larger than $0.2$', then you can be sure,
from the upper bound on the total variation distance, that none of your
HT-approximation predictions can be off by more than $10$~percentage
points from your true prediction.  This might be accurate enough for
your purposes.  Heuristically at least, you might conclude that
$\abs{\fstar_j - \ftilde_j} \ll 0.1$ for each $j$, on the grounds that
it would be unlikely for the entire error of the HT approximation to
concentrate on a single label.

Note the overtly subjective nature of the HT approximation, which
depends on your assessment of your $\beta$.  Exchangeability constrains
the PMF $\fbX$, but only to within an uncountably infinite set of
candidates (see section~\ref{sec:param} and the Appendix).  Your
particular $\fbX$ from this set might imply a very large $\beta$, as in
the example at the end of the previous section.  But, of course, if you
actually had a particular $\fbX$ you would not need an approximation.
So the HT serves those analysts, presumably the vast majority, who have
qualitative beliefs about $\bX$, which might well emerge into sharper
focus when tested with specific questions.  Contemplating $\bh$ and
asking ``How big is my $\beta$\/?'' is one such question.  To accept
that your $\beta$ is not larger than $0.2$ implicitly narrows your
$\fbX$ to a subset of all possible exchangeable PMFs, without requiring
you to be more specific.

We can also turn the HT around, to use the contrapositive.  If, knowing
the HT, you choose \emph{not} to use the HT approximation then you must
believe that at least one of the `add one' histograms adjacent to $\bh$
is much more probable, \textit{a priori}, than another.  So it turns out
in this case that you have strong beliefs about $\bX$: your rejection of
the HT approximation has revealed this to you, even if you did not know
it beforehand.

If you do not like $\bftilde$ as an approximation to your $\bfstar$,
then perhaps this is because $\bftilde$ is flatter than you would like.
It is certainly flatter than the `maximum likelihood' (ML) approximation 
\begin{displaymath}
  \fhat_j \ldef \frac{ h_j }{ n } \qquad j = 1, \dots, k.
\end{displaymath}
This is obvious from writing $\bftilde$ as the shrinkage estimator
\begin{displaymath}
  \bftilde = \frac{ n }{ n + k } \, \bfhat + \frac{ k }{ n + k } \, (k^{-1} \bone)
\end{displaymath}
where $\bone$ is the vector of $k$ ones.  It is straightforward to show
that 
\begin{equation}\label{eq:TVhat}
  \TV{ \bftilde - \bfhat } = \frac{ k }{ n + k } \, \TV{ \bfhat - (k^{-1} \bone) } .
\end{equation}
Therefore a necessary condition for preferring the more concentrated
$\bfhat$ to $\bftilde$ is that $n$ is not large relative to $k$,
otherwise there would be no appreciable difference between the two PMFs.
I will define an `under-powered' study as one where $n$ is not large
relative to $k$; say, for concreteness, that $n$ is not at least $9k$.
So, to put the condition differently, perhaps more starkly, you can only
find yourself in the position of preferring the ML approximation to the
HT approximation if you have performed an under-powered study.

Limited resources often constrain us to perform under-powered studies.
Suppose you find yourself, in such a study, preferring the ML
approximation to the HT approximation.  According to the established
norms of science \citep[see, e.g.,][ch.~3]{ziman00}, you must be careful
to ensure that your preference is rooted in an assessment of the
accuracy of the two approximations, rather than in the personal approval
that derives from presenting a more concentrated prediction.

One possibility for such an assessment is to modify the proof of
the HT to derive an upper bound on $\TV{ \bfstar - \bfhat }$.  But there
is a technical difficulty, which is that $h_j = 0$ would involve a
division by zero, and a trivial upper bound of $1$.  So consider a
slightly modified ML approximation 
\begin{displaymath}
  \fhtilde_j \ldef \frac{ h_j \vee 1 }{ n + v } \qquad j = 1, \dots, k
\end{displaymath}
where `$\vee$' denotes the pairwise maximum, and $v$ is the number of labels for
which $h_j = 0$.  Then it is straightforward to adapt the proof of the
HT to give $\TV{ \bfstar - \bfhtilde } \leq \half \gamma$, where 
\begin{displaymath}
  \gamma \ldef \frac{ \max_j\fbH^{n+1}(\bh \oplus j) \cdot (h_j + 1) / (h_j \vee 1) }{\min_{j'} \fbH^{n+1}(\bh \oplus j') \cdot (h_{j'} + 1) / (h_{j'} \vee 1) } .
\end{displaymath}
Unfortunately this bound does not lend itself to a simple assessment
procedure.  Unlike for $\beta$, it is no longer simply a case of
comparing histograms, which might be done semi-qualitatively, but,
instead, of attaching a probability to each histogram in order to make a
quantitative adjustment based on $\bh$.  However, a simple upper bound
on $\gamma$ can be found using the smallest and largest values of $\bh$,
denoted $\hlo$ and $\hhi$ respectively: 
\begin{displaymath}
  \gamma \leq \beta \, \frac { (\hlo + 1) / (\hlo \vee 1) }{ (\hhi + 1) / \hhi } .
\end{displaymath}
But since this upper bound is larger than $\beta$, and likely to be a
lot larger, it cannot serve the purpose of justifying the choice of the
ML approximation over the HT approximation.  So, to continue the stark
theme from the previous paragraph, it is hard to make the case that the
ML approximation is better than the HT approximation, in under-powered
studies where they differ.  From this I conclude that, in the absence of
strong beliefs for which $\beta$ is large, we should prefer the HT
approximation to the ML approximation.

\section{Population inference}
\label{sec:pop}

Usually, the sample of size $n$ is collected to make an inference about
the population of size $m$.  For example, the Mayor would like to know
the proportion of the electorate who are `very happy' with his
incumbency, measured on a Likert scale of $1$ (very unhappy) to $5$
(very happy).  Budgetary and time restrictions often limit the size of
the sample, such that the sample fraction $n / m$ is small.

\subsection{Point prediction}
\label{sec:point}

For a point prediction I will use the conditional expectation; i.e.\ the
Bayes estimate under quadratic loss \citep[][section~1.3, adopt a
similar approach to population prediction]{ghosh97}.  Denote the target
random quantity of interest as 
\begin{displaymath}
  P_j \ldef \frac{1}{m} \sum_{i=1}^m \indic_{X_i \deq x^{(j)}} \qquad j = 1, \dots, k ,
\end{displaymath}
the proportion of the population with label $j$.  Taking expectations,
\begin{equation}\label{eq:mayor}
\begin{split}
  \E^* (P_j)
& = \frac{1}{m} \left( h_j + \sum_{i=n+1}^m \Pr^*(X_i \deq x^{(j)}) \right) \\[0.25ex]
& = \left( \frac{n}{m} \right) \fhat_j + \left(1 - \frac{n}{m} \right)  \fstar_j \\[0.75ex]
& \approx \left( \frac{n}{m} \right) \fhat_j + \left(1 - \frac{n}{m} \right)  \ftilde_j \rdef \ptilde_j \qquad j = 1, \dots, k,
\end{split}
\end{equation}
using an asterisk to indicate conditioning on $\bX_{1:n} \deq
\bx_{1:n}$, as before.  The approximation in the final line uses the HT
to replace $\bfstar$ by $\bftilde$, on the condition that $\beta$ is
small.  This approximation can easily be adapted to predictions of
functions of $\bX$, for example in applications where $\curlyX$ is a set
of numbers or vectors.

Eq.~\eqref{eq:mayor} shows that the prediction $\bptilde$ is a weighted
average of the ML approximation and the HT approximation, where the
weight is the sample fraction $n / m$.  As was discussed in
section~\ref{sec:impl}, in under-powered studies it is possible that
$\bftilde \not\approx \bfhat$, and the case was made there that we
should prefer $\bftilde$ to $\bfhat$ as an approximation to $\bfstar$.
Thus, if the sample fraction is very small, then the approximate
prediction for small $\beta$ should be $\bptilde \approx \bftilde$ and
not $\bptilde \approx \bfhat$.

To be absolutely clear, the barchart displaying $\bptilde$ should be
titled `Prediction based on a survey of size $n$', and the vertical axis
should be labelled `Proportion of the electorate'.  By way of contrast,
the sample bar chart should be titled `Survey of size $n$', and the
vertical axis should be labelled `Number'.  This distinction should
alert other analysts to the different meanings of the two possibly
similar but definitely not identical-looking figures.

\subsection{Inference involving groups}
\label{sec:groups}

Often the population can be divided into distinct groups, where each
member belongs to exactly one group.  For example, the electorate can be
divided into men and women, or divided over the product of several
factors, such as sex, age, and ethnicity.  Another possibility is to
stratify the population according to the binned values of one or more
continuous values.

Where there are groups, it is natural to implement exchangeability
within groups but not across groups, as discussed in \citet{lindley81}.
This need not involve any additional statistical modelling under an
additional condition, given immediately below in~\eqref{eq:cind}.

Let there be $s$ groups indexed by $g$, so that the random quantity of interest is
\begin{displaymath}
  P_j = \sum_{g = 1}^s \frac{m_g}{m} P_{gj} \qquad j = 1, \dots, k,
\end{displaymath}
where $m_g$ is the size of group $g$ (taken as known, at least
approximately), and $P_{gj}$ is the proportion of group $g$ which has
label $j$.  The crucial simplifying assumption is that the samples are
sufficiently large, and the groups sufficiently distinct, that
\begin{equation}\label{eq:cind}
  \E^*(P_{gj}) \approx \E( P_{gj} \given \bX^g_{1:n_g} \deq \bx^g_{1:n_g}) \rdef \E^*_g(P_{gj}) \qquad
  \begin{cases}
    & g = 1, \dots, s \\ & j = 1, \dots, k .
  \end{cases}
\end{equation}
In other words, in the presence of the sample from group $g$, the
information from the other groups in the sample conveys effectively no
additional information about group $g$.  Then
\begin{equation}\label{eq:strat}
  \E^*(P_j) \approx \sum_g \frac{ m_g }{ m } \E^*_g(P_{gj}) 
  \approx \sum_g \frac{ m_g }{ m } \, \ptilde_{gj} \qquad j = 1, \dots, k,
\end{equation}
where $\ptilde_{gj}$ is the HT approximation from \eqref{eq:mayor},
applied just to group $g$.

If the approximation in \eqref{eq:cind} holds, then \eqref{eq:strat} is
an approximation to the Bayes estimate under quadratic loss,
because it is an approximation to the conditional expectation.  If
$\beta_g$ is small for each group then it is a good approximation,
according to the HT.  \citet{heiberger14} describe an attractive and
efficient way of displaying both the group predictions and the overall
population prediction in one figure.

But what about when \eqref{eq:cind} seems dubious, perhaps because the
sample size $n_g$ is small for some groups?  In this case
\eqref{eq:strat} is a point prediction, but not a good approximation to
an optimal one.  A better prediction could be constructed, but at the
cost of specifying a more restrictive statistical model for $\bX$ and
its group structure.  The standard approach would be a hierarchical
model such as 
\begin{equation}\label{eq:hier}
\begin{aligned}
  X^g_1, \dots, X^g_{m_g} \given \bphi, \theta & \iid f_{X \given \phi}(\cdott; \phi_g) \quad g = 1, \dots, s \\
  \phi_1, \dots, \phi_s \given \theta & \iid f_{\phi \given \theta}(\cdott; \theta) \\
  \theta & \sim f_\theta(\cdot)
\end{aligned}
\end{equation}
\citep[see][ch.~5]{gelman14}.  As required, this model is exchangeable
within each group, but not across groups.  It has the additional
parsimonious property that the group parameters are themselves
exchangeable, which allows learning across groups; i.e., `borrowing
strength' when $n_g$ is small.  This model requires three distributions
to be specified: the conditional distributions $f_{X \given \phi}$ and
$f_{\phi \given \theta}$, and the marginal distribution $f_\theta$, and
inference may require Monte Carlo methods, although these would be
standard \citep[see, e.g.,][]{lunn13}.

It is understandable that many analysts will prefer to accept the
approximate nature of \eqref{eq:cind} and \eqref{eq:strat}, than commit
to a specific set of distributional choices to implement
\eqref{eq:hier}, and the additional risk of a mathemetical or computing
error.  In particular, where sampling is not blocked by group, or where
response rates are low in some groups, analysts may prefer to increase
the sample sizes of under-sampled groups retrospectively, in order to
make \eqref{eq:cind} more appropriate.

\section{Parametric models}
\label{sec:param}

\citet{brewer09} identify prediction-based estimation with
\emph{parametric} modelling for the population.  Yet, as the previous
two sections showed, there need be no overt mention of parameters at all
in an exchangeable model for the population, without or with groups.
This section explores this issue in more detail.

Consider an IID parametric model of the form
\begin{displaymath}
  X_1, \dots, X_m \iid f_X(\cdot \,; \theta) \qquad \theta \in \Omega \subset \R^d .
\end{displaymath} 
This model, which is exchangeable for each $\theta$, is the dominant
statistical model of our time.  It asserts, first, that the proportion
of the population in the set $C \subset \curlyX$ can be modelled as 
\begin{displaymath}
  \sum_{j=1}^k \indic_{x^{(j)} \in C} \cdot f_X(x^{(j)}; \theta) ,
\end{displaymath}
for some value of $\theta$; second, that the selection process can be
modelled by random sampling with replacement (which is `ignorable').
The usual practice for population inference with this parametric model
would be to derive a point estimate of $\theta$ from the sample (e.g.\
using Maximum Likelihood), and plug this estimate into the model to
derive a point estimate of any population quantities of interest.
Variation in the estimate due to the small sample size can be assessed
using the bootstrap, or some other asymptotic approximation.

Now consider this IID parametric model in the light of the Finite
Representation Theorem (FRT) stated and proved in the Appendix.  The FRT
shows that there is a bijection between the set of exchangeable PMFs for
$\bX$ and the unit simplex
\begin{equation}\label{eq:simplex}
  \simplex^{c-1} \ldef \Set{ \bw \in \R^c : w_r \geq 0, \sum_{r=1}^c w_r = 1 } ,
\end{equation}
where $c$ is given in \eqref{eq:defc}, and represents the number of
different histograms that can be created with $m$ objects and $k$
labels.  We can think of $\bw \in \simplex^{c-1}$ as the universal
parameter of an exchangeable PMF.  It is straightforward to show that,
for an IID model 
\begin{displaymath}
  w_r = \text{M} \big( \bu^{(r)}; \bp(\theta) \big) \quad r = 1, \dots, c,
\end{displaymath}
where `$\text{M}$' denotes the Multinomial PMF, $\bu^{(r)}$ is the $r$th
histogram in some ordering, and $$\bp(\theta) \ldef \big( f_X(x^{(1)};
\theta), \dots, f_X(x^{(k)}; \theta) \big).$$  So for the IID model, the
set of exchangeable models does not occupy the whole of the universal
parameter space $\simplex^{c-1}$, but only a $d$-dimensional manifold
within $\simplex^{c-1}$.

This then is the general characteristic of `parametric' models for
exchangeable $\bX$: they can be represented by $\bw \in \Phi$, where
$\Phi$ is a strict subset of $\simplex^{c-1}$, and $\simplex^{c-1}
\setminus \Phi$ represents all of the exchangeable PMFs that are ruled
out \textit{a priori}.  It is reasonable in this context to denote
exchangeable models with no restrictions on $\simplex^{c-1}$ as
`non-parametric'.  Thus the previous sections showed that it is possible
to construct a non-parametric population model, and the presence of
parameters should not be conflated with population modelling, as
\citet{brewer09} have done.

The HT from section~\ref{sec:HT} is a non-parametric result: it holds
for all exchangeable $\bX$.  But consider what happens where a
parametric model is specified.  For example, consider the IID case where
$f_X$ is Normal and $\theta = (\mu, \sigma^2)$; in this case,
$\bp(\theta)$ is always bell-shaped.  When the analyst considers the
`add one' histograms adjacent to $\bh$, she will see some that are
heading towards bell-shaped, and some that are heading away from it.
Since she has ruled out non-bell-shaped $\bp$'s, she will attach more
probability to the more-bell-shaped than the less-bell-shaped add-one
histograms, i.e.\ her $\beta$ in the HT will be larger than zero.
Therefore, although the HT applies to all exchangeable sequences, the
analyst with a parametric model has a predisposition to rate some of the
`add one' histograms as more probable than others, and hence a
predisposition to believe that $\beta$ is not small, and that the HT
approximation $\bftilde$ might be a poor approximation to $\bfstar$.

Having said that, there is nothing to stop the analyst from comparing
her parametric prediction with the prediction she would make using the
HT approximation $\bftilde$.  The difference between these two
predictions will reveal the extent to which her beliefs have shaped her
prediction.  If she finds that the difference is large, and if she is
not confident in the beliefs that she has incorporated into her
parametric model, then she will need to reconsider.

There can be no harm in always requiring analysts who use parametric
models also to provide a prediction based on the HT approximation.  In
situations where the analyst has control over the number of labels, the
assessment of section~\ref{sec:impl} suggests setting, say, $k \assign
\lfloor n / 9 \rfloor$.  In this case there will be little difference
between the ML approximation and the HT approximation.  If the analyst
requires higher resolution in the labels, i.e.\ a larger $k$, then she
should accept that the price of an under-powered study is a flatter
prediction, and favour the HT approximation over the ML approximation.

\section{Merging the labels}
\label{sec:merge}

One immediate implication of the exchangeability of $\bX$ is that if
$Z_i$ is any function of $X_i$, then $\bZ \ldef (Z_1, \dots, Z_m)$ is
also exchangeable.  In particular, $Z_i$ could be a re-labelling of
$X_i$ in which two or more of the labels are merged; for example,
$x^{(1)}, x^{(2)}, x^{(3)}$ might no longer be distinguished, but all
labelled as $z^{(1)}$, with $z^{(j)} = x^{(j+2)}$ for the remaining
labels.

This presents an interesting conundrum, for the HT approximation.  If
$\bX$ is an exchangeable sequence then there are two routes to a
prediction for the merged labels: predict-then-sum, and
merge-then-predict.  The value is the same in both cases.  But for the
HT approximation, the value is different in the two cases.
Predict-then-sum adds $k$ to the total count, while merge-then-predict
adds less than $k$: $k-2$ in the example above.  Both routes are valid,
because if the original sequence is exchangeable, then the relabelled
sequence is exchangeable.  But, as the next illustration shows, the two
HT approximations can be completely different.

Consider a questionnaire in which $30$ questions are asked, and each one
is answered on a Likert scale of $1$ to $5$.  A sample of $n = 1000$ is
collected.  This is $n = 1000$ responses for $k = 5^{30} \approx
10^{21}$ distinct labels.  A prediction is required for one particular
question.  The HT approximation for predict-then-sum would increment the
count from $n = 10^3$ to $n + k \approx 10^{21}$.  The single question
margin is found by summing across all labels with the same
single-question label, and these margins would be almost completely
uniform.  For the other route, merge-then-predict would reduce the
number of labels to $5$ by merging the labels of the other questions.
Then the HT approximation would increment the count from $1000$ to
$1005$, and the prediction would be almost identical to the ML
approximation (see section~\ref{sec:impl}).

I think everyone would agree that the histogram-shaped HT-approximation
for merge-then-predict is better than the uniform HT-approximation for
predict-then-sum.  In the terms of the HT, we need to understand why
$\beta$ for predict-then-sum is so much larger than $\beta$ for
merge-then-predict.  Furthermore, this needs to be an \textit{a priori}
argument, since we make it without reference to any particular
histogram~$\bh$.

One obvious point, to be made and then passed over, is that for
predict-then-sum it is impossible to follow the procedure outlined at
the start of section~\ref{sec:impl}, because the number of add-one
histograms adjacent to $\bh$ is literally astronomical.  So in fact
$\beta$ could never be assessed for predict-then-sum.  However, an
\textit{a priori} argument in favour of merge-then-predict would obviate
the need to actually check the add-one histograms, and so the
impossibility of the procedure would not signify.

Reassuringly, there is an \textit{a priori} argument that the $\beta$
for predict-then-sum is large: larger than $2$, which is the point
at which the total variation bound is trivial.  I will use balls-in-bins
for clarity.  You arrange $1000$ balls across $10^{21}$ bins in some
fashion, representing the histogram $\bh$.  Now consider each of the
add-one histograms adjacent to $\bh$.  For a tiny fraction of these
histograms, you are adding a ball to a bin which already has a ball in
it, but for the rest you are adding a ball to an empty bin.  Under any
reasonably vague beliefs you would be very surprised indeed if, in your
$1001$-ball histogram over $10^{21}$ bins, the extra ball ended up in an
already-occupied bin, since the occupied bins are a tiny fraction of the
total bins.  Or, to put it differently, you would have to have
extraordinarily strong beliefs to attach roughly the same probability to
the extra ball ending up in an occupied bin, as an unoccupied one.
So, my claim is that, for reasonably vague beliefs
\begin{displaymath}
  \frac{ \Pr(\text{$1001$th ball in unoccupied bin}) }{ \Pr(\text{$1001$th ball in occupied bin}) }
\end{displaymath}
is larger than $3$ (much larger, I suggest), in which case $\beta > 2$,
and $\TV{\bfstar - \bftilde} = 1$, making the HT approximation useless
for prediction.

The universal practice in presenting questionnaire results is to display
the questions marginally; i.e.\ to merge the labels of the other
questions.  The only thing that I would change about this practice is to
display the HT approximation $\bftilde$ rather than the ML approximation
$\bfhat$, as described in section~\ref{sec:impl}.  Or, for inferences
with large sample fractions or with groups, to use the generalisations
given in section~\ref{sec:pop}.

\setcounter{equation}{0}
\renewcommand{\theequation}{A\arabic{equation}}
\section*{Appendix}

This is a self-contained Appendix deriving the Finite Representation
Theorem for exchangeable $\bX \ldef (X_1, \dots, X_m)$ where each $X_i$
has a finite realm $\curlyX \ldef \Set{ x^{(1)}, \dots, x^{(k)}}$.  It
provides precise statements about the probability mass functions (PMFs)
in \eqref{eq:fbXm} and \eqref{eq:fbXn}, and the form of $\fbH^{n+1}(\bh
\oplus j)$ in the Histogram Theorem.  It also clarifies the nature of
the parameter in an exchangeable PMF.

Let
\begin{displaymath}
  \curlyU \ldef \Set{ \bu^{(1)}, \dots, \bu^{(c)} }
\end{displaymath}
be the set of all possible histograms of $m$ objects over $k$ labels, where
\begin{equation}\label{eq:defc}
  c = \binom{m + k - 1}{k - 1}
\end{equation}
according to the elegant `stars and bars' construction of
\citet[sec.~II.5]{feller68}.  Assume, for concreteness, that these
histograms are in lexicographic order, indexed by $r$.

The definition of exchangeable $\bX$ was given in the main text, in
\eqref{eq:EXdef}.  Two equivalent formulations are that the probability
mass function (PMF) of $\bX$ is a symmetric function of $\bx$, and that
the PMF depends only on the histogram of $\bx$.  The key result is as
follows.

\begin{theorem}[Finite Representation Theorem, FRT]\label{thm:FRT}\samepage
  $f_{\bX_{1:n}}$ is the marginal PMF of an exchangeable $m$-sequence $\bX$ if
  and only if it has the form
\begin{equation}\label{eq:FRT}
  f_{\bX_{1:n}}(\bx_{1:n}) = \sum_{r=1}^c \frac{ \curlyH^n (\bh; \bu^{(r)}) }{ \curlyM_{\bh} } \cdot w_r
\end{equation}
where $\bh \ldef (h_1, \dots, h_k)$ and $h_j$ is the number of
$\bx_{1:n}$ with label $x^{(j)}$, $\curlyH^n$ is the multivariate
Hypergeometric distribution (see eq.~\ref{eq:curlyH} below),
$\curlyM_{\bh}$ is the Multinomial coefficient (see
eq.~\ref{eq:curlyM} below), and $\bw \ldef(w_1, \dots, w_c)$ lies in the
${(c-1)}$-dimensional unit simplex (see eq.~\ref{eq:simplex} in the main text).
\end{theorem}

The multivariate Hypergeometric distribution represents a random draw of
$n$ balls without replacement from an urn containing $u_j$ balls of each
of $k$ different colours, making $m$ balls altogether.  The probability
of drawing $h_j$ balls of colour $j$ for $j = 1, \dots, k$ is
\begin{equation}\label{eq:curlyH}
  \curlyH^n(\bh; \bu) \ldef \binom{m}{n}^{-1} \prod_{j=1}^k \binom{u_j}{h_j} \quad \text{if $\bh \leq \bu$ and $\sumnl_j h_j = n$}
\end{equation}
and zero otherwise.  The Multinomial coefficient
\begin{equation}\label{eq:curlyM}
  \curlyM_{\bh} \ldef \frac{ n! }{ h_1! \cdots h_k! }
\end{equation}
represents the number of distinct sequences $\bx_{1:n}$ which have the
same histogram $\bh$.

The FRT states that every marginal distribution of $\bX$ has the form of
a mixture over random draws without replacement from $m$-urns of
different compositions.  Setting $n \assign m$ shows that there is a
bijection between the set of exchangeable $\fbX$ and the
$(c-1)$-dimensional simplex indexed by $\bw$.  In this sense $\bw$ is
the universal parameter of an exchangeable $\bX$, and the
$(c-1)$-dimensional simplex is its parameter space.  From a Bayesian
point of view, $w_r$ is the prior probability that $\bX$ has histogram
$\bu^{(r)}$.

To prove the FRT, start by conditioning on $\bH$, the
histogram of $\bX$, and apply the Law of Total Probability:
\begin{equation}\label{eq:LTP}
  \fbX(\bx) = \sum_{r=1}^c f_{\bX \given \bH}(\bx \given \bu^{(r)}) \cdot \fbH^m (\bu^{(r)}) .
\end{equation}
The first term in the summand must be zero unless $\bh(\bx) =
\bu^{(r)}$, where $\bh(\cdot)$ is the histogram function.  For $\fbX$ to
be exchangeable,
\begin{equation}\label{eq:fX|H}
   \fbH^m(\bu^{(r)}) > 0 \implies f_{\bX \given \bH}(\bx \given \bu^{(r)}) = \indic_{\bh(\bx) \deq \bu^{(r)}} \, \big/ \, \curlyM_{\bu^{(r)}}  ,
\end{equation}
so that the probability is shared equally over all $\bx$ with the same
histogram $\bu^{(r)}$.  When $\fbH^m(\bu^{(r)}) = 0$, the value of
$f_{\bX \given \bH}(\bx \given \bu^{(r)})$ is immaterial in
\eqref{eq:LTP}, so we can take it to be \eqref{eq:fX|H} for all $r$.
Hence, for exchangeable $\bX$,
\begin{equation}\label{eq:fbXm1}
  \fbX(\bx) = \sum_{r=1}^c \frac{ \indic_{\bh(\bx) \deq \bu^{(r)}} }{ \curlyM_{\bu^{(r)}} } \cdot \fbH^m(\bu^{(r)}) = \fbH^m\{ \bh(\bx) \} \, \big/ \, \curlyM_{\bh(\bx)} .
\end{equation}
Eq.~\eqref{eq:fbXm1} is \eqref{eq:fbXm} in the main text.

Now to marginalise out $\bx_{(n+1):m}$, starting from
\begin{equation}\label{eq:fbXn1}
  f_{\bX_{1:n}}(\bx_{1:n}) 
= \sum_{\bx_{(n+1):m}} \fbX(\bx) 
= \sum_{r=1}^c \left\{ \sum_{\bx_{(n+1):m}} 
  \frac{ \indic_{\bh(\bx) \deq \bu^{(r)}} }{ \curlyM_{\bu^{(r)}} } \right\} \cdot \fbH^m(\bu^{(r)}) ,
\end{equation}
from the first equality in \eqref{eq:fbXm1}.  Let $\bh$ and $\bh'$
denote the histograms of $\bx_{1:n}$ and $\bx_{(n+1):m}$, respectively,
so that $\bh(\bx) = \bh + \bh'$.  Then the term in curly brackets
simplifies as
\begin{equation}
\begin{split}
\sum_{\bx_{(n+1):m}} 
  \frac{ \indic_{\bh(\bx) \deq \bu^{(r)}} }{ \curlyM_{\bu^{(r)}} }
& = \sum_{\bh'} \curlyM_{\bh'} \, \frac{ \indic_{\bh + \bh' \deq \bu^{(r)}} }{ \curlyM_{\bu^{(r)}} } \\
& = \sum_{\bh'} \frac{ \curlyM_{\bh'} }{ \curlyM_{\bu^{(r)}} } \, \indic_{\bh' \deq \bu^{(r)} - \bh} \\
& = \frac{ \curlyM_{\bu^{(r)} - \bh} }{ \curlyM_{\bu^{(r)}} } 
= \frac{ \curlyH^n (\bh; \bu^{(r)}) }{ \curlyM_{\bh} } \, ,
\end{split}
\end{equation}
after some re-arranging for the final equality.  Inserting this result
into \eqref{eq:fbXn1} completes the proof of the FRT.  In the statement
of the FRT I write $w_r \ldef \fbH^m(\bu^{(r)})$ to emphasise that any
point in the $(c-1)$-dimensional unit simplex is a valid choice for
$\bw$.

Rearranging \eqref{eq:FRT} gives
\begin{equation}\label{eq:fbXn2}
f_{\bX_{1:n}}(\bx_{1:n}) = \sum_{r=1}^c \curlyH^n (\bh; \bu^{(r)}) \cdot \fbH^m(\bu^{(r)}) \, \big/ \, \curlyM_{\bh} \rdef \fbH^n(\bh) \, \big/ \, \curlyM_{\bh} 
\end{equation}
which is \eqref{eq:fbXn} in the main text.  The critical probabilities
in the Histogram Theorem have the form
\begin{equation}
  \fbH^{n+1}(\bh \oplus j) = \sum_{r=1}^c \curlyH^{n+1} (\bh \oplus j; \bu^{(r)}) \cdot \fbH^m(\bu^{(r)}) ,
\end{equation}
i.e.\ they are deduced from $\fbH^m$ alone, albeit in a complicated way.

One of the subtle features of the FRT is that if $n < m$ then
$f_{\bX_{1:n}}$ in \eqref{eq:fbXn2} is a strict subset of the set of all
exchangeable PMFs for a sequence of length $n$.  This is because
$f_{\bX_{1:n}}$ must be `extendable' to an exchangeable $f_{\bX_{1:m}}$
\citep[see, e.g.,][]{aldous85}.  Thus the FRT is the finite analogue of
De~Finetti's Representation Theorem, which states that only a mixture of
IIDs is arbitrarily extendable.  The FRT states that only a mixture of
$m$-urns is extendable to $\bX_{1:m}$.  Heuristically, the FRT `proves'
De~Finetti's Representation Theorem, because letting $m$ increase for
fixed $n$ sends the multivariate Hypergeometric distribution towards the
Multinomial distribution, which represents sampling with replacement,
and therefore IID structure.  \citet{freedman77}, \citet{diaconis80} and
\citet[ch.~1]{schervish95} have more details, while \citet{heath76}
consider the binomial case, where $k = 2$.

\bibliography{statistics}

\begin{thebibliography}{}

\bibitem[Aldous, 1985]{aldous85}
Aldous, D. (1985).
\newblock Exchangeability and related topics.
\newblock In {\em Ecole d'Ete St Flour 1983}, pages 1--198. Springer Lecture
  Notes in Math.~1117.
\newblock Available at
  \url{http://www.stat.berkeley.edu/~aldous/Papers/me22.pdf}.

\bibitem[Brewer, 2002]{brewer02}
Brewer, K. (2002).
\newblock {\em Combined Survey Sampling Inference}.
\newblock Arnold, London, UK.

\bibitem[Brewer and Gregoire, 2009]{brewer09}
Brewer, K. and Gregoire, T. (2009).
\newblock Introduction to survey sampling.
\newblock {\em Sample Surveys: Design, Methods and Applications}, 29A:9--37.

\bibitem[Diaconis, 1977]{diaconis77}
Diaconis, P. (1977).
\newblock Finite forms of {de~Finetti's} theorem on exchangeability.
\newblock {\em Synthese}, 36(2):271--281.

\bibitem[Diaconis, 1988]{diaconis88}
Diaconis, P. (1988).
\newblock Recent progress on {de~Finetti's} notions of exchangeability.
\newblock In {\em Bayesian Statistics 3}, pages 111--125. Oxford University
  Press, Oxford, UK.
\newblock With discussion and rejoinder.

\bibitem[Diaconis and Freedman, 1980]{diaconis80}
Diaconis, P. and Freedman, D. (1980).
\newblock Finite exchangeable sequences.
\newblock {\em The Annals of Probability}, 8(4):745--764.

\bibitem[Edwards et~al., 1963]{edwards63}
Edwards, W., Lindman, H., and Savage, L. (1963).
\newblock Bayesian statistical inference for psychological research.
\newblock {\em Psychological Review}, 70(3):193--242.

\bibitem[Ericson, 1969]{ericson69}
Ericson, W. (1969).
\newblock Bayesian models in sampling finite populations.
\newblock {\em Journal of the Royal Statistical Society, Series~B},
  31(2):195--224.
\newblock With discussion, 224--232.

\bibitem[Feller, 1968]{feller68}
Feller, W. (1968).
\newblock {\em An Introduction to Probability Theory and its Applications}.
\newblock John Wiley \& Sons, New York NY, USA, 3rd edition.

\bibitem[Freedman, 1977]{freedman77}
Freedman, D. (1977).
\newblock A remark on the difference between sampling with and without
  replacement.
\newblock {\em Journal of the American Statistical Association}, 72:681.

\bibitem[Gelman et~al., 2014]{gelman14}
Gelman, A., Carlin, J., Stern, H., Dunson, D., Vehtari, A., and Rubin, D.
  (2014).
\newblock {\em Bayesian Data Analysis}.
\newblock Chapman and Hall/CRC, Boca Raton, FL, USA, 3rd edition.

\bibitem[Ghosh and Meeden, 1997]{ghosh97}
Ghosh, M. and Meeden, G. (1997).
\newblock {\em Bayesian Methods for Finite Population Sampling}.
\newblock Chapman \& Hall, London, UK.

\bibitem[Heath and Sudderth, 1976]{heath76}
Heath, D. and Sudderth, W. (1976).
\newblock {De Finetti's} theorem on exchangeable variables.
\newblock {\em The American Statistician}, 30(4):188--189.

\bibitem[Heiberger and Robbins, 2014]{heiberger14}
Heiberger, R. and Robbins, N. (2014).
\newblock Design of diverging stacked bar charts for {Likert} scales and other
  applications.
\newblock {\em Journal of Statistical Software}, 57(5).
\newblock Available online, \url{http://www.jstatsoft.org/v57/i05}.

\bibitem[Kingman, 1978]{kingman78}
Kingman, J. (1978).
\newblock Uses of exchangeability.
\newblock {\em The Annals of Probability}, 6(2):183--197.

\bibitem[Lindley and Novick, 1981]{lindley81}
Lindley, D. and Novick, M. (1981).
\newblock The role of exchangeability in inference.
\newblock {\em The Annals of Statistics}, 9(1):45--58.

\bibitem[Little, 2004]{little04}
Little, R. (2004).
\newblock To model or not to model? {Competing} modes of inference for finite
  population sampling.
\newblock {\em Journal of the American Statistical Association}, 99:546--556.

\bibitem[Lunn et~al., 2013]{lunn13}
Lunn, D., Jackson, C., Best, N., Thomas, A., and Spiegelhalter, D. (2013).
\newblock {\em The {BUGS} Book: A Practical introduction to Bayesian Analysis}.
\newblock CRC Press, Boca Raton FL, USA.

\bibitem[Murphy, 2012]{murphy12}
Murphy, K. (2012).
\newblock {\em Machine Learning: A Probabilistic Perspective}.
\newblock MIT Press, Cambridge MA, USA.

\bibitem[Schervish, 1995]{schervish95}
Schervish, M. (1995).
\newblock {\em Theory of Statistics}.
\newblock Springer, New York NY, USA.
\newblock Corrected 2nd printing, 1997.

\bibitem[Ziman, 2000]{ziman00}
Ziman, J. (2000).
\newblock {\em Real Science: What it is, and what it means}.
\newblock Cambridge, UK: Cambridge University Press.

\end{thebibliography}

\end{document}